\documentclass{elsart}

\usepackage{amsfonts,amsmath,amssymb}

\begin{document}


\begin{frontmatter}

\title{Numerical analysis of a nonlocal parabolic problem resulting
from thermistor problem\corauthref{corA}}

\author{Moulay Rchid Sidi Ammi}
\ead{sidiammi@mat.ua.pt} \quad
\author{Delfim F. M. Torres}
\ead{delfim@ua.pt} \corauth[corA]{Research Report CM06/I-25,
Dep. Mathematics, Univ. Aveiro, July 2006.}

\address{Department of Mathematics,
      University of Aveiro,
      3810-193 Aveiro, Portugal}

\date{}

\begin{abstract}
We analyze  the spatially semidiscrete piecewise linear finite
element method for a nonlocal parabolic equation resulting from
thermistor problem. Our approach is based on the properties of the
elliptic projection defined by the bilinear form associated with the
variational formulation of the finite element method. We assume
minimal regularity of the exact solution that yields optimal order
error estimate. The full discrete backward Euler method and the
Crank-Nicolson-Galerkin scheme are also considered. Finally, a
simple algorithm for solving the fully discrete problem is proposed.
\end{abstract}


\end{frontmatter}


\section{Introduction}

We study the numerical approximation by the finite element scheme of
the nonlinear problem
\begin{equation}
\label{11}
\begin{gathered}
\frac{\partial u}{\partial t}-\nabla . (k(u)\nabla u) = \lambda
\frac{f(u)}{ \big ( \int_{\Omega} f(u)\, dx \big )^{2}} ,
\mbox{ in } \Omega \times ]0;\overline{t}[, \\
u = 0 \quad \mbox{ on }  \partial \Omega \times ]0;\overline{t}[,  \\
u(0)= u_0  \quad \mbox{ in }  \Omega,
\end{gathered}
\end{equation}
where $\Omega$ is a bounded domain in $\mathbb{R}^{2}$,
$\overline{t}$ is a positive fixed real, $f$ and $k$ are functions
from $\mathbb{R}$ to $\mathbb{R}$ satisfying the hypotheses
$(H1)-(H2)$ below, $\lambda$ is a positive parameter and $\nabla$
denotes the gradient with respect to the $x-$variables. The time
evolution model \eqref{11} describes the temperature profile of a
thermistor device with electrical resistivity $f$, see
\cite{ac,bl,ci1,es1,es3,es2,lac1,lac2,psx,tza}; the dimensionless
parameter $\lambda$ can be identified with the square of the applied
potential difference $V$ at the ends of the conductor. The system
\eqref{11} has been the subject of a variety of investigations in
the last decade. Existence of weak solutions to problems related
with the thermistor problem is proved in \cite{jr}, where the
mathematical treatment of this system apparently appears for the
first time. In \cite{es} the problem \eqref{11} for the special case
$k=1$ is considered, and then a backward Euler
time-semidiscretization method for the approximation of its solution
is proposed and analyzed. In this paper we propose a finite element
method to construct numerical approximations of the solutions of
problem \eqref{11} for the case when $k$ is different from the
identity. The formulation of the finite element method is standard
and it is based on a variational formulation of the continuous
problem. There is a vast literature on finite element methods for
nonlinear elliptic and parabolic problems. For example, we mention
the work \cite{dew} on the porous media equations, which are similar
to the Joule heating problem \cite{alm}. Compared to a standard
semilinear equation, the main challenge here is the nonstandard
nonlocal nonlinearity on the right-hand side of the partial
differential equation \eqref{11}.

On the other hand, error bounds are normally expressed in terms of
norms of the exact solution of the problem. It is well known that
the required regularity of the exact solution can be attained by
assuming enough regularity of data, sometimes supplemented with
compatibility conditions, see \cite{al,el,yx}. We then use
sufficient conditions in terms of the data of the problem and its
solution $u$ that yield error estimates (see hypotheses (H1)-(H3)
below).


\section{Main results and organization of the paper}

We denote by $(\cdot, \cdot)$ and $\|\cdot\|$ respectively the inner
product and the norm in $L^{2}= L^{2}(\Omega)$, by $\|\cdot\|_{s}$
the norm in the Sobolev space $H^{s}(\Omega)$, by $c$ some generic
positive constant which may depend upon the data and whose value may
vary from step to step. In Section~\ref{sec:sdp} we study spatially
semidiscrete approximations of \eqref{11} by the finite element
method. The approximate solution is sought in the piecewise linear
finite element space
$$
S_{h}=S_{h}(\Omega)= \left\{ \chi \in C(\Omega): \chi /_{e} \,
\,linear, \forall \, e \in T_{h}; \chi /_{\partial \Omega } \,  =0
\right\},
$$
where $\{ T_{h}\}_{h}$ is a family of regular triangulations of
$\Omega$, with $h$ denoting the maximum diameter of the triangles of
$ T_{h}$. As a model for our analysis we first consider the
corresponding semidiscrete Galerkin finite element method, which
consists in finding $u_{h}(t)\in S_{h}$ such that
\begin{equation}
\label{13}
\begin{gathered}
\left(u_{h,t}, \chi\right) + \left(k(u_{h})\nabla u_{h}, \nabla
\chi\right)
= \frac{\lambda}{ \big ( \int_{\Omega} f(u_{h})\, dx \big )^{2}}(f(u_{h}), \chi), \\
u_{h}(0)= u_{0h}  \, ,
\end{gathered}
\end{equation}
$\forall \chi \in S_{h}$, $t \in J=(0,\overline{t})$, and where
$u_{0h} \in S_{h}$ is a given approximation of $u_{0}$. Similar
discretization techniques have been analyzed for various linear and
nonlinear evolution problems (cf. e.g. \cite{el}). This method
\eqref{13} may be written as a system of ordinary differential
equations. In fact, let $\{ \phi \}_{j=1}^{N_{h}}$ be the standard
pyramid basis of $S_{h}$. Write $u_{h}(x, t)= \sum_{j=1}^{N_{h}}
\alpha_{j}(t)\phi_{j}(x)$, where $(\alpha_{j})_{1\leq j \leq N_{h}}$
are the real coefficients to be determined.  Then, \eqref{13} can be
written as
$$
A \alpha'(t)+ B(\alpha)\alpha(t)= \widetilde{f}(\alpha), \quad t \in
J \, , \quad \alpha(0)= \gamma,
$$
where $\gamma$ is the vector of nodal values of $u_{0h}$,
$\widetilde{f}(\alpha)= (\widetilde{f}_{1}(\alpha), \ldots,
\widetilde{f}_{N_{h}}(\alpha))^{T}$ with
$$
\widetilde{f}_{j}(\alpha)= \frac{\lambda}{ \big ( \int_{\Omega}
f(\sum_{l=1}^{N_{h}} \alpha_{l}(t)\phi_{l})\, dx \big
)^{2}}\left(f\left(\sum_{l=1}^{N_{h}} \alpha_{l}(t)\phi_{j}\right),
\phi_{k}\right),
$$
and $A= (a_{jk})_{1\leq j, k \leq N_{h}}$ and
$B(\alpha)=(b_{jk}(\alpha))_{1\leq j, k \leq N_{h}}$ are,
respectively, the corresponding mass and stiffness matrices:
$$
a_{jk}=\left( \phi_{j}, \phi_{k} \right) \, , \quad b_{jk}(\alpha) =
\left(k\left(\sum_{l=1}^{N_{h}} \alpha_{l}(t)\phi_{l}\right)\nabla
\phi_{j}, \nabla \phi_{k}\right) \, .
$$

We shall assume the following general assumptions on the given data:

\begin{description}

\item[(H1)] $f: \mathbb{R} \rightarrow \mathbb{R}$
is a locally Lipschitzian function  and $f(u) \geq \sigma > 0$ for
all $u \in \mathbb{R}$;

\item[(H2)] $k$ is a twice derivable function verifying: there exist
positive constants
 $k_{1}, k_{2}$ and $c$ such that $0<k_{1}\leq k(u)
\leq k_{2}$, $|k'(u), \,k''(u)| \leq c$;

\item[(H3)] $u \in L^{\infty}(0, T, H^{2}(\Omega)\bigcap W^{1, \infty}(\Omega))$ and $u_{0} \in H^{2}(\Omega)$.

\end{description}
It is shown in \cite[Theorem~2.1]{ci2} that the regularity
assumption $(H3)$ is satisfied if the data is smooth and compatible.
The matrix $A$ is always definite positive. Further, hypotheses (H1)
and (H2) assure that the matrix $B(\alpha)$ is also positive
definite. Assumption (H3) is used in order
 to prove Lemma~\ref{lem11} and then to show a $o(h^{2})$ error
estimate. It is useful to introduce the interpolation operator
$I_{h}: C(\Omega) \rightarrow S_{h}$ defined by
$$
I_{h}v= \sum_{j=1}^{N_{h}}v(P_{j})\phi_{j}(x),
$$
where $\{ P_{j}\}_{j=1}^{N_{h}}$ are the interior vertices of
$T_{h}$. For completeness, we assume the following standard
interpolation error estimates : for $v \in H^{2}(\Omega)\cap
H_{0}^{1}(\Omega)$ there exists a positive constant $c>0$ such that
$$
\|I_{h}v-v\| \leq ch^{2}\|v\|_{2} \mbox{ and }
 \|\nabla (I_{h}v-v)\| \leq ch \|v\|_{2},
$$
holds. We also assume the property \cite{vt}: for some integer
$r\geq 2$ and small $h$,
\begin{equation}
\label{aa} \inf_{\chi \in S_{h}} \{ \| v-\chi\|+h \|
\nabla(v-\chi)\|  \} \leq ch^{s} \|v \|_{s}, \mbox{ for } 1\leq s
\leq r,
\end{equation}
when $v \in H^{s}\bigcap H_{0}^{1}$.
 We finally suppose that the family of triangulations is such that the
inverse estimate \cite{vt}
\begin{equation}
 \| \nabla \chi \| \leq
ch^{-1}\| \chi \| \quad \forall \chi \in S_{h}
\end{equation}
is satisfied. In the existing literature the error estimates for the
finite element method are normally expressed in terms of norms of
the exact solution of the problem and are usually derived for
solutions that are sufficiently smooth (cf. e.g. \cite{el}). To
estimate the error in the semidiscrete problem \eqref{13} we split
the error as $u_{h}-u= (u_{h}-\widetilde{u}_{h})+
(\widetilde{u}_{h}-u) = \theta + \rho$, where $\widetilde{u}_{h}$
denotes the standard elliptic projection in $S_{h}$ of the exact
solution $u$ defined by
\begin{equation}
\label{24} \left(k(u(t))\nabla (\widetilde{u}_{h}-u), \nabla \chi
\right)=0, \quad \forall \chi \in S_{h}, \quad t\geq 0.
\end{equation}
It is well known (see \cite{esd,vt}) that if the regularity
hypothesis (H3) for $u$ holds, then $\widetilde{u}_{h}$ has the
following approximation properties:
\begin{lem}
\label{lem11} Under the regularity hypotheses (H1)-(H3), one has
$$
\|\rho(t)\|+h\|\nabla \rho(t)\| \leq c(u) h^{2},
$$
$$
\|\rho_{t}(t)\|+h\|\nabla \rho_{t}(t)\| \leq c(u) h^{2},
$$
where $c(u)$ is a constant independent of $t \in J$.
\end{lem}
\begin{lem}
\label{lem12} Let $\widetilde{u}_{h}$ be defined by \eqref{24}.
Then, $\|\nabla \widetilde{u}_{h}(t)\|_{L^{\infty}} \leq c(u)$,
$t\in J$.
\end{lem}

For the finite element method \eqref{13}, if $u_{0h}$ is chosen such
that $\|u_{0h}-u_{0}\| \leq ch^{2}\|u_{0}\|_{2}$, we prove in
Section~\ref{sec:sdp} an error estimate of the form
$\|u_{h}(t)-u(t)\| \leq c(u)h^{2}$ (Theorem~\ref{thm21}). The
corresponding estimate for the gradient is also proved
(Theorem~\ref{thm22}). In Section~\ref{sec:cdc} we show that our
approach for the semidiscrete Galerkin finite element method also
applies to fully discrete schemes. We consider the backward Euler
method for the discretization in time of \eqref{11}: letting $\tau$
to be the time step, $U^{n}$ the approximation in $S_{h}$ of $u(t)$
at $t=t_{n}=n\tau$, i.e. $U^{n}= \sum_{j=1}^{N_{h}}
\alpha_{j}^{n}\phi_{j}$, where $(\alpha_{j}^{n})_{1 \leq j \leq
N_{h}}$ are the unknown real coefficients,  $\partial_{n} U^{n} =
\frac{U^{n}-U^{n-1}}{\tau }$, $n= 0, 1, 2, \ldots,$ the numerical
method is defined by
\begin{equation}
\label{111}
\begin{gathered}
( \partial_{n} U^{n}, \chi ) +(k(U^{n})\nabla U^{n}, \nabla \chi) =
\frac{\lambda}{ \big ( \int_{\Omega} f(U^{n})\, dx \big
)^{2}}(f(U^{n}), \chi) ,
\quad \forall \chi \in S_{h}, \\
U^{0}= u_{0h}  .
\end{gathered}
\end{equation}
For this scheme we prove (Theorem \ref{thm31}), under the same
regularity requirements as in Section~\ref{sec:sdp}, that
$$
\|U^{n}- u(t_{n})\| \leq c(u) (h^{2}+\tau),
 \mbox{ for } t_{n} \in J=(0, \overline{t}) .
$$
In order to obtain higher accuracy in time, in
Section~\ref{sec:GNGs} we investigate an alternative way to obtain
an $o(h^{2}+\tau^{2})$ error bound using the basic
Crank-Nicolson-Galerkin scheme: applying the techniques of
sections~\ref{sec:sdp} and \ref{sec:cdc} we prove
(Theorem~\ref{thm}) an error estimate of the form
$$
\|U^{n}- u(t_{n})\| \leq c(u) (h^{2}+\tau^{2}), \mbox{ for } t_{n}
\in J=(0, \overline{t}).
$$
Finally, in Section~\ref{sec:alg} we propose a simple algorithm for
solving the fully discrete problem.


\section{Semidiscrete problem}
\label{sec:sdp}

In this section we obtain an error estimate for $u$ and the
associated estimate for the gradient. The proofs use a splitting of
the error based on the elliptic projection $\widetilde{u}_{h}$
\eqref{24}. We may define the semidiscrete problem on a finite
interval $J=(0, \overline{t}]$ of time.

\begin{thm}
\label{thm21} Let $u$ and $u_{h}$ be the solutions of \eqref{11} and
\eqref{13}, respectively. Then, under the hypotheses (H1)-(H3), we
have: $\|u_{h}(t)- u(t)\| \leq  c(u) h^{2}$, for  $t \in J$,
provided that $\|u_{0h}- u_{0}\| \leq  c h^{2}$.
\end{thm}

\begin{pf}
Owing to Lemma~\ref{lem11} and to the decomposition of the error as
sum of two terms $u_{h}-u= (u_{h}- \widetilde{u}_{h})
+(\widetilde{u}_{h}-u)=\theta+ \rho$, it suffices to treat $\theta =
u_{h}- \widetilde{u}_{h}$. We have from the equations satisfied by
$u_{h}$ and $\widetilde{u}_{h}$ that
\begin{align*}
\begin{split}
 &(\theta_{t}, \chi)+\left(k(u_{h})\nabla \theta, \nabla
\chi\right) \\
& =(u_{h, t}, \chi)+\left(k(u_{h})\nabla u_{h}, \nabla \chi\right)
-\left(\widetilde{u}_{h,t}, \chi\right)
-\left(k(u_{h})\nabla \widetilde{u}_{h}, \nabla \chi\right)\\
& =\frac{\lambda}{(\int f(u_{h}))^{2}} (f(u_{h}), \chi )- (\rho_{t},
\chi )- (u_{t}, \chi) - \left(k(u)\nabla \widetilde{u}_{h}, \nabla
\chi\right)
+ \left(( k(u)-k(u_{h}))\nabla \widetilde{u}_{h}, \nabla \chi\right)\\
&  =\frac{\lambda}{(\int f(u_{h}))^{2}} (f(u_{h}), \chi )
-(\rho_{t}, \chi ) - (u_{t}, \chi) - (k(u)\nabla {u}, \nabla \chi )
+ \left(( k(u)-k(u_{h}))\nabla \widetilde{u}_{h}, \nabla \chi\right)\\
&=\frac{\lambda}{(\int f(u_{h}))^{2}} (f(u_{h}), \chi ) -
\frac{\lambda}{(\int f(u))^{2}} (f(u), \chi) + \left((
k(u)-k(u_{h}))\nabla \widetilde{u}_{h}, \nabla \chi\right)
-(\rho_{t}, \chi )\\
& =\frac{\lambda}{(\int f(u_{h}))^{2}} (f(u_{h})-f(u)), \chi ) +
\left(\frac{\lambda}{(\int f(u_{h}))^{2}} - \frac{\lambda}{(\int
f(u))^{2}}\right) (f(u), \chi ))\\
 & + \left((
k(u)-k(u_{h}))\nabla \widetilde{u}_{h}, \nabla \chi\right) -
(\rho_{t},\chi )\\
 &=\frac{\lambda}{(\int f(u_{h}))^{2}}
(f(u_{h})-f(u)), \chi ) + \left(( k(u)-k(u_{h}))\nabla
\widetilde{u}_{h}, \nabla \chi\right) -
(\rho_{t},\chi )\\
& + \frac{\lambda}{(\int f(u_{h}))^{2} (\int f(u))^{2}}
\left(\int_{\Omega} (f(u_{h})-f(u))dx \right ) \left( \int_{\Omega}
(f(u_{h})+f(u)) dx \right )
 (f(u), \chi ).
\end{split}
 \end{align*}
Thus, setting $\chi = \theta$, using the hypotheses (H1)-(H3),
Lemma~\ref{lem11}, and Young's inequality,
\begin{equation*}
\begin{split}
\frac{1}{2}\frac{d}{dt}\|\theta\|^{2} +k_{1}\| \nabla \theta\|^{2}
&\leq c (\|u_{h}-u\|( \|\theta\|+ \|\nabla \theta\|)+ \|\rho_{t}\|
\|\theta\|) \\
&\leq k_{1} \|\nabla \theta\|^{2}+ c
(\|\theta\|^{2}+\|\rho\|^{2}+\|\rho_{t}\|^{2}).
\end{split}
\end{equation*}
By integration, we get $\|\theta (t)\|^{2} \leq \|\theta (0)\|^{2} +
c \int_{0}^{t}(\|\theta \|^{2} +\|\rho\|^{2}+\|\rho_{t}\|^{2}) ds$
and it follows by Gronwall's Lemma that
$$
\|\theta (t)\|^{2} \leq c\|\theta (0)\|^{2}+
c\int_{0}^{t}(\|\rho\|^{2}+\|\rho_{t}\|^{2}) ds
$$
or
\begin{equation}
\label{211} \|\theta (0)\| \leq \|u_{0h}-u_{0}\|+
\|\widetilde{u}_{h}(0)-u_{0}\|\leq \|u_{0h}-u_{0}\|+
ch^{2}\|u_{0}\|_{2}.
\end{equation}
We then get the desired conclusion:
$$
\|\theta (t)\| \leq c\|u_{0h}-u_{0}\|+ c(u)h^{2} \leq c(u)h^{2}.
$$
\end{pf}

We now derive for the standard Galerkin method, from the weak
formulation of the parabolic problem and using the inverse property,
the optimal order error estimate for the gradient.

\begin{thm}
\label{thm22} Let $u$ and $u_{h}$ be, respectively, the solutions of
\eqref{11} and \eqref{13}. Under hypotheses (H1)-(H3), if $u_{0h}$
is chosen such that $\|u_{0h}-u_{0}\| \leq ch^{2} \|u_{0}\|_{2}$,
then
$$
\|\nabla u_{h}(t)- \nabla u(t)\| \leq  c(u) h, \mbox{ for } t \in J.
$$
\end{thm}

\begin{pf}
We have
\begin{align}
\label{a1} \|\nabla u_{h}(t) - \nabla u(t)\| & \leq \|\nabla
(u_{h}(t)- \chi)\|
+ \|\nabla \chi- \nabla u(t)\| \nonumber \\
& \leq ch^{-1} \| u_{h}(t)-  \chi \|+\|\nabla \chi- \nabla u(t)\|\\
&\leq ch^{-1} \| u_{h}(t)-  u(t)\|+ ch^{-1}( \|\chi - u(t)\|
+h\|\nabla \chi - \nabla u(t)\|). \nonumber
\end{align}
By the approximation assumption \eqref{aa} we know that, with
suitable $\chi \in S_{h}$,
$$
\|\chi - u(t)\|+h \|\nabla \chi - \nabla u(t)\| \leq
ch^{2}\|u(t)\|_{2}.
$$
Then, using \eqref{a1}, we get:
$$
\|\nabla u_{h}(t)- \nabla u(t)\| \leq ch^{-1}\| u_{h}(t)-  u(t)\|
+ch\|u(t)\|_{2}.
$$
Theorem~\ref{thm21} yields the intended conclusion: $\|\nabla
u_{h}(t)- \nabla u(t)\| \leq c(u)h$.
\end{pf}


\section{The completely discrete case}
\label{sec:cdc}

We now turn our attention to the fully discrete scheme based on the
backward Euler method: find $U^{n} \in S_{h}$ such that \eqref{111}
holds. Existence result of \eqref{111} is a simple consequence of
Brouwer's fixed point theorem. Here we obtain an error estimate for
the scheme. For $h$, $\tau$ and $\frac{\tau}{h}$ small enough, we
prove uniqueness.

\begin{thm}
\label{thm31} Let $u$ and $U^{n}$ be solutions of \eqref{11} and
\eqref{111} respectively, with $u_{0h}$ chosen such that
$\|u_{0}-u_{0h}\| \leq ch^{2}$. Under the required regularity
(H1)-(H3), there exists a constant $c$ such that, for $t_{n} \in J$
and small $\tau$, we have
$$
\|U^{n}-u(t_{n})\| \leq c(u) (h^{2}+\tau) .
$$
Moreover, for sufficiently small $\tau$, $h$ and $\frac{\tau }{h}$,
there is a unique solution $U^{n}$ of the complete discrete scheme
\eqref{111}.
\end{thm}

\begin{pf}
We use the partitioning of the error
\begin{equation}
\label{ep} U^{n}-u^{n}= (U^{n}-\widetilde{U}^{n})+
(\widetilde{U}^{n}-u^{n}) = \theta_{n}+ \rho_{n} ,
\end{equation}
with $u^{n}=u(t_{n})$, $\widetilde{U}^{n}=
\widetilde{u}_{h}(t_{n})$, where $\widetilde{u}_{h}$ is the elliptic
projection of $u^{n}$ defined by \eqref{24}. By virtue of
Lemma~\ref{lem11}, it suffices to bound $\theta_{n}$. We have for
$\chi \in S_{h}$ that
\begin{align*}
( \partial_{n}\theta_{n} , \chi ) & +(k(U^{n})\nabla \theta_{n},
\nabla \chi)\\
& =  ( \partial_{n} U^{n}, \chi ) +(k(U^{n})\nabla U^{n}, \nabla
\chi) -  ( \partial_{n} \widetilde{U}^{n}, \chi )
-(k(U^{n})\nabla \widetilde{U}^{n}, \nabla \chi) \\
& = \frac{\lambda}{ \big ( \int_{\Omega} f(U^{n})\, dx \big
)^{2}}(f(U^{n}), \chi) -(u_{t}^{n},\chi) -( \partial_{n}
\widetilde{U}^{n}- u_{t}^{n}, \chi) \\
& \qquad -(k(u^{n})\nabla \widetilde{U}^{n}, \nabla
\chi)-((k(U^{n})-k(u^{n}))\nabla \widetilde{U}^{n}, \nabla \chi)
\end{align*}
and, in view of the elliptic projection \eqref{24} and the equation
of the continuous problem \eqref{11}, we can write:
 \begin{align*}
( \partial_{n}\theta_{n} , \chi ) &+ (k(U^{n})\nabla \theta_{n},
\nabla \chi)\\
& = \frac{\lambda}{ \big ( \int_{\Omega} f(U^{n})\, dx \big
)^{2}}(f(U^{n}), \chi)-\frac{\lambda}{ \big ( \int_{\Omega}
f(u^{n})\, dx \big )^{2}}(f(u^{n}), \chi)\\
& \qquad -( \partial_{n} \rho_{n}, \chi )  - ( \partial_{n} {u}^{n}-
u_{t}^{n}, \chi)- ((k(U^{n})-k(u^{n}))\nabla \widetilde{U}^{n},
\nabla \chi)\\
& = \frac{\lambda}{ \big ( \int_{\Omega} f(U^{n})\, dx \big
)^{2}}(f(U^{n})-f(u^{n}), \chi)\\
& \qquad + ( \frac{\lambda}{ \big ( \int_{\Omega} f(U^{n})\, dx \big
)^{2}}-\frac{\lambda}{ \big (
\int_{\Omega} f(u^{n})\, dx \big )^{2}})(f(u^{n}), \chi)\\
& \qquad -( \partial_{n} \rho_{n}, \chi )-( \partial_{n} {u}^{n} -
u_{t}^{n}, \chi)- ((k(U^{n})-k(u^{n}))\nabla \widetilde{U}^{n},
\nabla \chi).
\end{align*}
Choosing $\chi = \theta_{n}$ and using the fact that $\nabla
\widetilde{U}^{n}$, $u^{n}$, $U^{n}$ are bounded, we get:
\begin{multline*}
\frac{1}{2} \partial_{n}\| \theta_{n}\|^{2}+ k_{1} \|\nabla
\theta_{n}\|^{2}\\
\leq c \|U^{n}-u^{n}\|(\|\theta_{n}\|+\|\nabla \theta_{n}\|) +
(\|\partial_{n} \rho_{n}\|+\|\partial_{n} {u}^{n}- u_{t}^{n}\| )\|
\theta_{n}\|.
\end{multline*}
Hence, by Young's inequality,
\begin{equation}
\label{112}
\begin{gathered}
\partial_{n}\| \theta_{n}\|^{2}+ k_{1} \|\nabla \theta_{n}\|^{2}
\leq c(\| \theta_{n}\|^{2}+ \| \rho_{n}\|^{2}+\|
\partial_{n} \rho_{n}\|^{2}+\|\partial_{n} {u}^{n}-
u_{t}^{n}\|^{2}).
\end{gathered}
\end{equation}
Introducing the notation $R_{n}= \| \rho_{n}\|^{2}+\|
\partial_{n} \rho_{n}\|^{2}+\|\partial_{n} {u}^{n}-
u_{t}^{n}\|^{2}$ we write \eqref{112} in the form $\partial_{n}\|
\theta_{n}\|^{2}+ k_{1} \|\nabla \theta_{n}\|^{2} \leq c(\|
\theta_{n}\|^{2}+R_{n})$, and it results that
$$
(1-c\tau )\| \theta_{n}\|^{2} \leq \| \theta_{n-1}\|^{2}+c\tau
R_{n}.
$$
 For  $\tau < \frac{1}{c}$, we
have $\| \theta_{n}\|^{2} \leq \frac{1}{1-c\tau}\|
\theta_{n-1}\|^{2}+\frac{c\tau R_{n}}{1-c \tau}$. Using the fact
that  $\frac{1}{1-c\tau} \thickapprox 1 +c \tau$ for $\tau$
sufficiently small, it follows that
 $\| \theta_{n}\|^{2} \leq (1+c\tau
)\| \theta_{n-1}\|^{2}+c\tau R_{n}$. By induction, we get:
\begin{equation}
\label{113}
\begin{gathered}
\| \theta_{n}\|^{2} \leq (1+c\tau )^{n}\|
\theta_{0}\|^{2}+c\tau \sum_{j=1}^{n}(1+c\tau )^{n-j}R_{j} \\
\leq c \| \theta_{0}\|^{2}+ c\tau \sum_{j=1}^{n} R_{j}, \mbox{ for }
t_{n} \in J.
\end{gathered}
\end{equation}
We now recall that from Lemma~\ref{lem11}
$$
\|\rho_{j} \| \leq c(u) h^{2} \, , \quad \|\partial_{n} \rho_{j} \|=
\left\|\frac{1}{\tau } \int_{t_{j-1}}^{t_{j}}\rho_{t} ds\right\|
\leq c(u) h^{2}.
$$
On the other hand, we have
$$
\left\|\partial_{n} {u}^{j}-u_{t}^{j}\right\|^{2} = \left\|
\frac{1}{\tau } \int_{t_{j-1}}^{t_{j}}(s-t_{j-1}) u_{tt}(s) ds
\right\| \leq c(u) \tau .
$$
It yields that $R_{j} \leq c(u) (h^{2}+\tau )^{2}$. Taking the above
estimates together with \eqref{113} and \eqref{211} we prove the
intended bound for $\theta_{n}$:
$$\| \theta_{n}\| \leq c
\|u_{0h}-u_{0}\|+ c(u)(h^{2}+\tau) \leq  c(u)(h^{2}+\tau).$$

It remains to prove the second part of the theorem (uniqueness). Let
$U^{n}=X$ and $U^{n}=Y$ be two solutions of the fully discrete
problem:
\begin{align*}
(X &- Y, \chi )+\tau \left(k(X)\nabla X-k(Y)\nabla Y, \nabla \chi\right)\\
&= \frac{\lambda \tau }{(\int f(X))^{2}}(f(X)-f(Y),
\chi)+\left(\frac{\lambda \tau }{(\int f(X))^{2}}-\frac{\lambda
\tau}{(\int f(Y))^{2}}\right)(f(Y), \chi).
\end{align*}
Taking $\chi = X-Y$, we have
\begin{align*}
& \|X-Y\|^{2}+\tau (k(X)\nabla (X-Y), \nabla (X-Y))= \frac{\lambda
\tau }{(\int f(X))^{2}}(f(X)-f(Y), X-Y)\\
&+\left(\frac{\lambda \tau }{(\int f(X))^{2}}-\frac{\lambda \tau
}{(\int f(Y))^{2}}\right)(f(Y), X-Y) -\tau \left((k(X)-k(Y))\nabla
Y, \nabla (X-Y)\right).
\end{align*}
Thus,
$$
\|X-Y\|^{2}+ \frac{1}{2}k_{1}\tau \|\nabla (X-Y)\|^{2}\leq
c\|X-Y\|^{2}(\tau +\tau \|\nabla Y\|_{L^{\infty}}^{2}).
$$
According to Lemma~\ref{lem12} we have
$$ \|\nabla Y\|_{L^{\infty}} \leq \|\nabla
\widetilde{u}_{h}\|_{L^{\infty}}+\|\nabla \theta_{n}\|_{L^{\infty}}
\leq c+ch^{-1}\|\nabla \theta_{n}\|.$$ Taking into account  the
estimate for $R_{n}$, we get
$$
k_{1}\|\nabla \theta_{n}\|^{2}\leq c\|\theta_{n-1}\|^{2}+\tau
\|\theta_{n}\|^{2}+\tau R_{n} \leq c(h^{2}+\tau )^{2},
$$
and we deduce that
$$
\tau \|\nabla Y\|_{L^{\infty}}^{2} \leq c \left(\tau +h^{2} +
\left(\frac{\tau}{h}\right)^{2}\right).
$$
Then, for sufficient small $\tau$, $h$ and $\frac{\tau }{h}$, we get
the uniqueness of the solution of the complete discrete scheme.
\end{pf}


\section{The Crank-Nicolson-Galerkin scheme}
\label{sec:GNGs}

This section is devoted to the study of the following
Crank-Nicolson-Galerkin scheme:
\begin{equation}\label{cng}
\begin{gathered}
\left( \partial_{n} U^{n}, \chi \right) +
\left(k(\overline{U}^{n})\nabla \overline{U}^{n}, \nabla \chi\right)
= \frac{\lambda}{ \big ( \int_{\Omega} f(\overline{U}^{n})\, dx \big
)^{2}}(f(\overline{U}^{n}), \chi), \,
\forall \chi \in S_{h}, \, t_{n}\in J, \\
U^{0}= u_{0h},
\end{gathered}
\end{equation}
with $\overline{U}^{n}=\frac{U^{n}+U^{n-1}}{2}$. Before proceeding
to the main result of the section -- an $o(h^{2} + \tau^{2})$ error
bound -- we need an auxiliary estimate.

\begin{lem}
\label{lem41} Besides the hypotheses (H1)-(H3), let us further
assume
\begin{description}
\item[(H4)] $u_{tt} \in L^{\infty}(0, T, H^{1}(\Omega))$.
\end{description}
Then, $\|\nabla \widetilde{u}_{h,tt}\| \leq c(u)$, where
$\widetilde{u}_{h}$ is the elliptic projection defined by
\eqref{24}.
\end{lem}

\begin{pf}
Differentiating \eqref{24} twice in time we find that
$$
(k(u)\nabla \widetilde{u}_{h,tt},\nabla \chi ) = (k(u)\nabla u_{tt},
\nabla \chi )- 2(k'(u) \nabla \rho_{t}, \nabla \chi )-( k''(u)
\nabla \rho, \nabla \chi ).
$$
Taking $\chi = \widetilde{u}_{h,tt}$, it follows that
$$k_{1}
\|\nabla \widetilde{u}_{h,tt}\|^{2} \leq c(u)(\|\nabla u_{tt}\| +
\|\nabla \rho_{t}\|+\|\nabla \rho\|) \|\nabla
\widetilde{u}_{h,tt}\|$$ and, recalling Lemma~\ref{lem11} and (H4),
we obtain the intended conclusion.
\end{pf}

\begin{thm}
\label{thm} Let $u$ and $U^{n}$ be respectively the solutions of
\eqref{11} and \eqref{cng}. Then, under the hypotheses (H1)-(H4), we
have $\| U^{n} -u(t_{n})\| \leq c(u)(h^{2}+\tau^{2})$, $t_{n}\in J$.
\end{thm}

\begin{pf}
Splitting the error by means of the elliptic projection as in
\eqref{ep}, one gets
\begin{align*}
\begin{split}
(\partial_{n} & \theta^{n}, \chi)
+ \left(k(\overline{U}^{n})\nabla \overline{\theta}^{n}, \nabla \chi\right) \\
& = \left( \partial_{n} U^{n}, \chi \right) +
\left(k(\overline{U}^{n})\nabla \overline{U}^{n}, \nabla \chi\right)
- \left( \partial_{n} \widetilde{U}^{n}, \chi \right)
- \left(k(\overline{U}^{n})\nabla \overline{\widetilde{U}}^{n}, \nabla \chi\right) \\
& = \frac{\lambda}{\left(\int_{\Omega} f(\overline{U}^{n})\, dx
\right)^{2}} \left(f(\overline{U}^{n}), \chi\right)-
\left(u_{t}^{n-\frac{1}{2}}, \chi\right)
- \left(\partial_{n} \widetilde{U}^{n}-u_{t}^{n-\frac{1}{2}}, \chi\right) \\
& \qquad - \left(k\left(u^{n-\frac{1}{2}}\right)\nabla
\widetilde{U}^{n-\frac{1}{2}}, \nabla \chi\right)
-\left(k(\overline{U}^{n})\nabla
\overline{\widetilde{U}}^{n}-k\left(u^{n-\frac{1}{2}}\right)
\nabla \widetilde{U}^{n-\frac{1}{2}}, \nabla \chi\right) \\
& =\frac{\lambda}{\left(\int_{\Omega}
f\left(\overline{U}^{n}\right)\, dx \right)^{2}}
\left(f(\overline{U}^{n}), \chi\right)
-\frac{\lambda}{\left(\int_{\Omega} f(u^{n-\frac{1}{2}})\, dx
\right)^{2}}
\left(f(u^{n-\frac{1}{2}}), \chi\right) \\
& \qquad - \left( \partial_{n}
\widetilde{U}^{n}-u_{t}^{n-\frac{1}{2}}, \chi\right) -
\big(\left(k(\overline{U}^{n})- k(u^{n-\frac{1}{2}})\right)
\nabla \overline{\widetilde{U}}^{n} \\
& \qquad  + k(u^{n-\frac{1}{2}})\nabla
\left(\overline{\widetilde{U}}^{n} -
\widetilde{U}^{n-\frac{1}{2}}\right), \nabla \chi\big).
\end{split}
\end{align*}

Setting $\chi = \overline{\theta}^{n}$, it follows from
$\left(\partial_{n} \theta^{n}, \overline{\theta}^{n}\right) =
\frac{1}{2}\overline{\partial}\|\theta^{n}\|^{2}$ and (H1)-(H3) that
\begin{align*}
\begin{split}
& \frac{1}{2}\overline{\partial} \|\theta^{n}\|^{2} + \mu \|\nabla
\overline{\theta}^{n} \|^{2} \leq  \frac{\lambda}{ \left(
\int_{\Omega} f(u^{n-\frac{1}{2}})\, dx \right)^{2}}|
(f(u^{n-\frac{1}{2}})-f(\overline{U}^{n}), \nabla
\overline{\theta}^{n} )| \\
& \qquad + \left|\frac{\lambda}{\big( \int_{\Omega}
f(u^{n-\frac{1}{2}})\, dx \big)^{2}}- \frac{\lambda}{\big (
\int_{\Omega} f(\overline{U}^{n})\, dx
\big)^{2}}(f(\overline{U}^{n}), \nabla
\overline{\theta}^{n})\right| \\
& \qquad + c \big( \|\partial_{n}
\widetilde{U}^{n}-u_{t}^{n-\frac{1}{2}}\| +
\|\overline{U}^{n}-u^{n-\frac{1}{2}}\|+\|\nabla
(\overline{\widetilde{U}}^{n}
-\widetilde{U}^{n-\frac{1}{2}})\|\big)\|\nabla
\overline{\theta}^{n}\| \\
&  \leq c \big( \|\partial_{n}
\widetilde{U}^{n}-u_{t}^{n-\frac{1}{2}}\| +
\|\overline{U}^{n}-u^{n-\frac{1}{2}}\|+\|\nabla
(\overline{\widetilde{U}}^{n} -\widetilde{U}^{n-\frac{1}{2}})\|
\big)\|\nabla \overline{\theta}^{n}\| .
\end{split}
 \end{align*}
Young's inequality gives
\begin{equation}
\label{43} \overline{\partial}\|\theta^{n}\|^{2} \leq c \big(
\|\partial_{n} \widetilde{U}^{n}-u_{t}^{n-\frac{1}{2}}\|^{2} +
\|\overline{U}^{n}-u^{n-\frac{1}{2}}\|^{2}+\|\nabla
(\overline{\widetilde{U}}^{n} -\widetilde{U}^{n-\frac{1}{2}})\|^{2}
\big).
\end{equation}
Estimating each term on the right hand side of the inequality
\eqref{43} separately, we have
\begin{equation}
\label{44} \|\overline{U}^{n}-u^{n-\frac{1}{2}}\| \leq
\|\overline{\theta}^{n}\|+\|\overline{\rho}^{n}\| +
\|\overline{u}^{n}- u^{n-\frac{1}{2}}\| \leq
\|\overline{\theta}^{n}\| + c(u) (h^{2}+\tau^{2}),
\end{equation}
\begin{equation}
\label{45} \|\partial_{n}
\widetilde{U}^{n}-u_{t}^{n-\frac{1}{2}}\|\leq \|\partial_{n}
\rho^{n}\| +\|\partial_{n} u^{n}-u_{t}^{n-\frac{1}{2}}\| \leq c(u)
(h^{2}+\tau^{2}),
\end{equation}
and by Lemma~\ref{lem41}
\begin{equation}
\label{46} \left\|\nabla (\overline{\widetilde{U}}^{n}
-\widetilde{U}^{n-\frac{1}{2}})\right\|  \leq c \tau
\int_{t_{n-1}}^{t_{n}} \|\nabla \widetilde{u}_{h, tt}\| \, ds \leq
c(u) \tau^{2}.
\end{equation}
Inequalities \eqref{43}--\eqref{46} together show that
\begin{align*}
\begin{split}
\overline{\partial}\|\theta^{n}\|^{2} &\leq c
\|\overline{\theta}^{n}\|^{2}+c(u)(h^{2}+ \tau^{2})^{2} \leq
\frac{c}{4}\| \theta^{n}+ \theta^{n-1}\|^{2}
+c(u)(h^{2}+ \tau^{2})^{2} \\
& \leq c \| \theta^{n}\|^{2}+ c \| \theta^{n-1}\|^{2}+c(u)(h^{2}+
\tau^{2})^{2} ,
\end{split}
\end{align*}
which gives $(1-c \tau) \| \theta^{n}\|^{2} \leq (1+c \tau) \|
\theta^{n-1}\|^{2} +c(u)\tau (h^{2}+ \tau^{2})^{2}$. Then, by
induction, we have that for small $\tau$
\begin{equation*}
\| \theta^{n}\| \leq c \| \theta^{0}\| + c(u)(h^{2}+ \tau^{2}) \leq
c \|u_{0h}-u_{0}\|+ c(u)(h^{2}+ \tau^{2}) \mbox{ for } t_{n} \in J .
\end{equation*}
Or $\|u_{0h}-u_{0}\| \leq c h^{2}$, and then $\| \theta^{n}\| \leq
c(u)(h^{2}+ \tau^{2})$, which leads, in view of Lemma~\ref{lem11},
to the intended result.
\end{pf}


\section{Method of solution}
\label{sec:alg}

We divide the interval $\Omega=(-1, 1)$ into $N$ equal finite
elements: $x_{0}=-1 < x_{1} < \ldots < x_{N}=1$. Let $(x_{j},
x_{j+1})$ be a partition of $\Omega$ and $x_{j+1}-x_{j}=h=
\frac{1}{N}$ the step length. By $S$ we denote a basis of the usual
pyramid functions:
\begin{equation*}
v_{j}(x)=
\begin{cases}
\frac{1}{h}x+(1-j)
 & \mbox{ on } [x_{j-1}, x_{j}],\\
- \frac{1}{h}x+(1+j)
 & \mbox{ on } [x_{j}, x_{j+1}],\\
 0 & otherwise.
\end{cases}
\end{equation*}
We first write \eqref{11} in variational form with $k=1$. We then
multiply \eqref{11} by $v_{j}$ (for $j$ fixed). We have, using the
boundary conditions, that
\begin{equation*}
\int_{\Omega} \frac{\partial u}{\partial t} v_{j} \, dx +
\int_{\Omega}  \nabla u \nabla v_{j} \, dx  = \lambda
\frac{\int_{\Omega}f(u) v_{j} \, dx}{ \big ( \int_{\Omega} f(u)\, dx
\big )^{2}}.
\end{equation*}
By the Crank-Nicolson approach we obtain
\begin{align}\label{eq17}
&\int_{\Omega} u^{n+1} v_{j} \, dx + \tau \int_{\Omega} \nabla
u^{n+1} \nabla v_{j} \, dx = \int_{\Omega} u^{n} v_{j} +\lambda \tau
\frac{\int_{\Omega}f(u^{n}) v_{j}\, dx}{ \big ( \int_{\Omega}
f(u^{n})\, dx \big )^{2}},
\end{align}
so the approximation $u^{n+1}$ to the function $u$ can be written as
$$u^{n+1}= \sum_{i=-1}^{N}\alpha_{i}^{n+1} v_{i},$$ where the
$\alpha_{i}^{n+1}$ are unknown real coefficients to be determined.
From \eqref{eq17} it is easy to obtain the following system of
$(N-1)$ linear algebraic equations:
\begin{multline}
\label{eq18}
  \left (\frac{h}{6} -\frac{\tau}{h} \right )\alpha_{j-1}^{n+1}
  +\left ( \frac{2h}{3}
+\frac{\tau}{2h} \right) \alpha_{j}^{n+1}
+ \left ( \frac{h}{6} -\frac{\tau}{h}\right) \alpha_{j+1}^{n+1}  \\
 =\frac{h}{6}\alpha_{j-1}^{n} + \frac{2h}{3}\alpha_{j}^{n}+
 \frac{h}{6}\alpha_{j+1}^{n}+
\frac{\lambda \tau \int_{\Omega}f(u^{n}) v_{j} \, dx}{ \big (
\int_{\Omega} f(u^{n})\, dx \big )^{2}}.
\end{multline}
Using the Dirichlet  boundary conditions, we have
\begin{equation}\label{eq19}
\begin{split}
\alpha_{-1}^{n+1}& = \alpha_{1}^{n+1} -(1+h)\alpha_{0}^{n+1},\\
\alpha_{-1}^{n} &= \alpha_{1}^{n} -(1+h)\alpha_{0}^{n},
\end{split}
\end{equation}
and
\begin{equation}\label{eq20}
\begin{split}
\alpha_{N}^{n+1}& = c \alpha_{N-1}^{n+1}, \\
\alpha_{N}^{n} &= c \alpha_{N-1}^{n} ,
\end{split}
\end{equation}
where $c = \frac{Nh}{(N-1)h-1}$. Substituting the expressions
\eqref{eq19} and \eqref{eq20} into \eqref{eq18}, we obtain the
system of equations
\begin{equation}\label{eq21}
\begin{split}
\left ( b-(1+h) \right )\alpha_{0}^{n+1} &+2a
\alpha_{1}^{n+1}\\
& =\left ( \frac{2h}{3}-(1+h) \right )\alpha_{0}^{n}+
\frac{h}{3}\alpha_{1}^{n}+\frac{\lambda \tau \int_{\Omega}f(u^{n})
v_{0} \, dx}{ \big ( \int_{\Omega} f(u^{n})\, dx \big )^{2}}, \, j=0\\
a\alpha_{j-1}^{n+1} +b \alpha_{j}^{n+1}
&+ a \alpha_{j+1}^{n+1} \\
 & = \frac{h}{6}\alpha_{j-1}^{n} + \frac{2h}{3}\alpha_{j}^{n}+
  \frac{h}{6}\alpha_{j+1}^{n}+
   \frac{\lambda \tau \int_{\Omega}f(u^{n}) v_{j} \, dx}{ \big (
\int_{\Omega} f(u^{n})\, dx \big )^{2}}, \, j=1 \ldots N-2, \\
a \alpha_{N-2}^{n+1} &+ (ac +b) \alpha_{N-1}^{n+1}\\
 &= \frac{h}{6}\alpha_{N-2}^{n}+\frac{h}{6}(c+4)\alpha_{N-1}^{n} +
\frac{\lambda \tau \int_{\Omega}f(u^{n}) v_{N-1} \, dx}{ \big (
\int_{\Omega} f(u^{n})\, dx \big )^{2}}, \, j= N-1,
\end{split}
\end{equation}
where $$ a= \frac{h}{6}- \frac{\tau}{h} \, , \quad  b=\frac{2h}{3}+
\frac{\tau}{2h},$$ and thus a simple algorithm for solving the fully
discrete problem.


\begin{ack}
The authors were supported by FCT (\emph{The Portuguese Foundation
for Science and Technology}): Sidi Ammi through the fellowship
SFRH/BPD/20934/2004, Torres through the R\&D unit CEOC. We are
grateful to two anonymous referees for useful remarks and comments.
\end{ack}



\end{document}